\documentclass[preprint,12pt]{elsarticle}

\usepackage{amssymb}
\usepackage{amscd,amsmath,amsfonts,graphicx}
\usepackage{amsthm}

\begin{document}
\baselineskip=22pt \centerline{\large \bf Absolute moments in terms
of characteristic functions} \vspace{1cm} \centerline{Gwo Dong Lin
and Chin-Yuan Hu} \centerline{Academia Sinica and National Changhua
University of Education} \vspace{1cm} \noindent {\bf Abstract.} The
absolute moments of probability distributions are much more
complicated than conventional ones. By using a direct and simpler
approach, we retreat P.\,L.\,Hsu's (1951, J.\,\,Chinese
Math.\,\,Soc., {Vol.\,\,1}, pp.\,\,257--280) formulas in terms of
the characteristic function (which have been ignored in the
literature) and provide
 some new results as well. The case of nonnegative random variables  is also investigated
through both characteristic function and Laplace--Stieltjes
transform. Besides, we prove that the distribution of a nonnegative
random variable with a finite fractional  moment can be completely
determined by a {\it proper} subset
 of the translated fractional moments.
 This improves significantly P.\,Hall's (1983, Z.\,W., Vol.\,62, 355--359) result for distributions on the right-half line.
\vspace{0.1cm}\\
\hrule
{\bf Short title:  Formulas of absolute moments}\\
{\bf Postal addresses:} Gwo Dong Lin, Institute of Statistical
Science, Academia Sinica, Taipei 11529, Taiwan. (E-mail:
 gdlin@stat.sinica.edu.tw)\\
  Chin-Yuan Hu, Department of Business
Education, National Changhua University of Education, Changhua
50058, Taiwan. (E-mail: buhuua@gmail.com)\\
\newpage

\noindent{\bf 1. Introduction}
\newcommand{\bin}[2]{
   \left(
     \begin{array}{@{}c@{}}
         #1  \\  #2
     \end{array}
   \right)          }

The absolute moments of probability distributions  play important
roles in both theoretical and applied fields (see, e.g., von Bahr
1965, Hall 1983,  Nabeya 1951, 1952, Barndorff-Nielsen and Stelzer
2005, and  Matsui and Pawlas 2016). Our main purpose in this paper
is to investigate the presentation of (fractional) absolute moments
of distributions in terms of their Fourier--Stieltjes or
Laplace--Stieltjes transforms.

We review first some important properties of the characteristic
function (ch.f.). Consider a random variable $X$ with distribution
$F(x)=\Pr(X\le x)$
 on ${\mathbb R}\equiv(-\infty,\infty),$ denoted $X\sim F,$  and let
$\phi$ be the ch.f.\,\,of $F,$ namely, the Fourier--Stieltjes
transform
\begin{eqnarray}
\phi(t)=\int_{-\infty}^{\infty}e^{itx}\,dF(x)=
\int_{-\infty}^{\infty}\cos(tx)\,dF(x)+i\int_{-\infty}^{\infty}\sin(tx)\,dF(x),\
\ t\in{\mathbb R},
\end{eqnarray}
or, $\phi(t)={\bf E}[\exp(itX)] ={\bf E}[\cos(tX)]+i\,{\bf
E}[\sin(tX)]\equiv \hbox{Re}(\phi(t))+i\,\hbox{Im}(\phi(t)),\
t\in{\mathbb R}.$ It is known that for any distribution $F,$ its
ch.f.\,\,$\phi$ in (1) always exists and uniquely determines $F$ by
the famous Weierstrass Approximation Theorem. Therefore, we can
recover  $F$ and derive its properties via the ch.f.\,\,$\phi$
theoretically.  In fact, we have the classical L$\grave{\hbox{e}}$vy
inverse formula:
\begin{eqnarray*}
F(x+h)-F(x)=\lim_{T\to\infty}\frac{1}{2\pi}\int_{-T}^T\frac{1-e^{-ith}}{it}e^{-itx}\phi(t)\,dt,
\end{eqnarray*}
provided that both $x$ and $x+h$ (with $h>0$) are {\it continuity
points} of $F$
 (see, e.g., Lukacs 1970, Chapter 3).
 Another explicit  expression is the following (improper) integral:
 for any {\it continuity point} $x$ of $F$ on ${\mathbb R},$
\begin{eqnarray}
F(x)=\frac{1}{2}-\frac{1}{\pi}\int_0^{\infty}{\hbox{Im}}(\phi(t)e^{-itx})\,\frac{dt}{t}=
\frac{1}{2}-\frac{1}{\pi}\lim_{T\to\infty}\int_0^{T}{\hbox{Im}}(\phi(t)e^{-itx})\,\frac{dt}{t}
\end{eqnarray}
(see Zolotarev 1957, or Kawata 1972, p.\,130, or Rossberg et
al.\,\,1985, p.\,45).

On the other hand, if $X\ge 0,$ then for any {\it continuity point}
$x$ of $F$ on ${\mathbb R}_+\equiv[0,\infty),$ we have the inverse
formulas via (2):
\begin{eqnarray}
F(x)=\frac{2}{\pi}\int_0^{\infty}(\sin xt){{\hbox{Re}}(\phi(t))}\,\frac{dt}{t},\ \ x\in{\mathbb R}_+,\\
F(x)=1-\frac{2}{\pi}\int_0^{\infty}(\cos xt){{\hbox{Im}}(\phi(t))}\,\frac{dt}{t},\ \ x\in{\mathbb R}_+
\end{eqnarray}
 (Laue 1983, 1986). This means that  for any nonnegative random variable,
 one of the real and imaginary parts of its ch.f.\,\,is enough to recover the corresponding distribution.

It is also well known that if the random variable $X\sim F$ on
${\mathbb R}$  has finite $n$th  moment $m_n={\bf E}[X^n]$ for some
positive integer $n,$ then its ch.f.\,\,$\phi$ is $n$-time
differentiable and ${\bf E}[X^n]=(-i)^n\phi^{(n)}(0)$ (in terms of
the ch.f.). Unlike the conventional moment, the absolute moment
$\mu_s={\bf E}[|X|^s]$ is, however, much more complicated. For
example, if $\mu_s$ is finite for an odd positive integer $s,$ then
we can derive
\begin{eqnarray*}\mu_s=\frac{1}{2\pi
i^{s+1}}\int_{-\infty}^{\infty}[\phi^{(s)}(t)-\phi^{(s)}(-t)]\,\frac{dt}{t}
\end{eqnarray*}
 (which corrects a result in Lukacs 1970,
p.\,26). In particular,
\begin{eqnarray}
\mu_1=-\frac{1}{\pi}\int_{-\infty}^{\infty}\hbox{Re}(\phi^{\prime}(t))\,\frac{dt}{t}.
\end{eqnarray}
Moreover, $\mu_1$ has another formula
$$\mu_1=\frac{2}{\pi}\int_0^{\infty}\frac{1-\hbox{Re}(\phi(t))}{t^2}\,dt,$$
which is a special case of Theorem 1  or Theorem 1$^{\prime}$ below
 (see also Remark 1 for  discussions on (5) and the absolute moment
 $\mu_{-1}$ of
negative order). It is worth mentioning that the distribution of $X$
with finite $\mu_s$ is completely determined by the set
\begin{eqnarray}M_s=\{{\bf E}[|X+a|^s]: a\in{\mathbb R}\}
\end{eqnarray}  of translated
absolute moments, where
 $s>0$ is not an even integer (Hall 1983; see also Riesz 1938, Rudin 1976 and Mattner 1992 for high-dimension case with different proofs).

We will focus on the formulas of (fractional) absolute moments of
positive {\it real} orders. The main results are stated in the next
section, and all the proofs are given in Section 4. The case of
nonnegative random variables  is also studied through
Laplace--Stieltjes transform. Besides, we improve significantly
Hall's (1983) characterization result (6) for distributions on the
right-half line (Theorems 7 and 8). In Section 3, we give the
necessary lemmas to prove the main results. Finally, in
 Section 5 some illustrative remarks are provided.

\newpage
\noindent{\bf 2.  Main results}

Let $X\sim F$ on ${\mathbb R}$ with ch.f.\,\,$\phi$ and define the
polynomial
\begin{eqnarray*}
P_n(t)=\sum_{k=0}^n\frac{1}{k!}\phi^{(k)}(0)\,t^k,\ \ \ t\in{\mathbb R},
\end{eqnarray*}
provided $\phi$ is $n$-time differentiable.  Hsu (1951) might be the
first one to derive formulas of absolute moments. For example, he
obtained the following results in terms of the ch.f.\,\,(by using
the orthogonality property of the Hermite polynomials), in which
$\lfloor{s}\rfloor$ denotes the largest integer less than or equal
to $s.$

\noindent {\bf Theorem 1} (Hsu 1951, Theorems 2.1 and 4.1). Let
$X\sim F$ on ${\mathbb
R}.$\\
(a) If $X$ has finite absolute moment $\mu_{2n+1}$ of positive odd
order ${2n+1},$ then
\begin{eqnarray}
\mu_{2n+1}=(-1)^{n+1}\frac{2}{\pi}(2n+1)!\int_0^{\infty}\frac{\hbox{Re}(\phi(t)-P_{2n}(t))}{t^{2n+2}}\,dt.
\end{eqnarray}
(b) If the fractional absolute moment $\mu_{n+r}<\infty$ for some
integer $n\ge 0$ and real $r\in(0,1),$
\begin{eqnarray}\mu_{n+r}=(-1)^{\lfloor{n/2}\rfloor
+1}\frac{2}{\pi}\Gamma(n+r+1)\sin\left(\left(\frac{n}{2}-\lfloor{\frac n2}\rfloor+\frac r2\right)\pi\right)
\int_0^{\infty}\frac{\hbox{Re}(\phi(t)-P_n(t))}{t^{n+r+1}}\,dt.\vspace{0.4cm}
\end{eqnarray}

\vspace{0.1cm}Note that (7) can be rewritten as
$$\mu_{2n+1}=(-1)^{n+1}\frac{2}{\pi}(2n+1)!\int_0^{\infty}t^{-2n-2}\biggl[\hbox{Re}(\phi(t))-
\sum_{k=0}^{n}\frac{(-1)^k}{(2k)!}m_{2k}t^{2k}\biggr]\,dt.$$
In order to get a better insight into this result, let us recall the following:\\
(i) for $n\ge 0,$ the even function
\begin{eqnarray*}
g_n(t)=(-1)^{n+1}\biggl(\cos t-\sum_{k=0}^n\frac{(-1)^k}{(2k)!}t^{2k}\biggr)\ge 0,\ \ \ t\in{\mathbb R},
\end{eqnarray*}
and if $n\ge 1,$ $g_n$ is increasing and convex on ${\mathbb R}_+;$\\
(ii) for $n\ge 1,$ the odd function
\begin{eqnarray*}
h_n(t)=(-1)^{n}\biggl(\sin
t-\sum_{k=1}^n\frac{(-1)^{k+1}}{(2k-1)!}t^{2k-1}\biggr)\ge 0,\ \ \
t\ge 0,
\end{eqnarray*}
increases on ${\mathbb R}_+.$ (See, e.g., Hu and Lin 2005, Lemma 4.)

 Moreover, we define the auxiliary functions:
\begin{eqnarray*}
G_n(t)
&=&(-1)^{n+1}\biggl(\hbox{Re}(\phi(t))-\sum_{k=0}^n\frac{(-1)^k}{(2k)!}m_{2k}t^{2k}\biggr),\ \ \
t\in{\mathbb R},\ \ n\ge 0,\\
H_n(t)
&=&(-1)^{n}\biggl(\hbox{Im}(\phi(t))-\sum_{k=1}^n\frac{(-1)^{k+1}}{(2k-1)!}m_{2k-1}t^{2k-1}\biggr),\ \ \
t\in{\mathbb R},\ \ n\ge 1,
\end{eqnarray*}
provided the moments are finite.  Also, define $h_0(t)=\sin t$ and
$H_0(t)={\bf E}[\sin (tX)],\ t\in{\mathbb R}.$

Combining (7) and (8), we have the following slight extension of
Theorem 1. Note that we don't assume the finiteness of expectations
in (9) because of the fact $g_n\ge 0,$  and hence both sides could
be infinite simultaneously.

\noindent {\bf Theorem 1$^{\prime}$.} Let $X\sim F$ on ${\mathbb
R}$, $n\ge 0$ be an integer and $2n<s<2n+2.$ Then
\begin{eqnarray}
\mu_{s}={\bf E}[|X|^s]=\frac{2}{\pi}\Gamma(s+1)\sin((s-2n)\pi/2)
\int_0^{\infty}\frac{{\bf E}[g_n(tX)]}{t^{s+1}}\,dt\ \ \ (\hbox{finite or infinite}).
\end{eqnarray}
If, in addition, $\mu_s<\infty$ for some $s\in(2n, 2n+2)$ with
integer $n\ge 0,$ then
\begin{eqnarray}
\mu_{s}=\frac{2}{\pi}\Gamma(s+1)\sin((s-2n)\pi/2)\int_0^{\infty}\frac{G_n(t)}{t^{s+1}}\,dt.
\end{eqnarray}

Theorems 1 and 1$^{\prime}$ and some of their special cases have
been rediscovered in the literature (see, e.g., von Bahr 1965, von
Bahr and Esseen 1965, Brown 1970, 1972,  Kawata 1972, Section 11.4,
Chung 2001, Section 6.2, and Chow and Teicher 1997, Section 8.4).
Surprisingly,  few of the above authors mentioned the original
article Hsu (1951). A direct and probably simpler proof of Theorems
1$^{\prime}$ will be given later.

Similarly, we have the following result for distributions on the
right-half line, which modifies a result of
 Brown (1970, Lemma 1) (noting that the function
$h_0(t)=\sin t,\ t\ge 0,$ is not of constant sign) and retreats
those of Laue (1986) in a simpler way. The latter applied the
calculus of fractional derivatives and integrals.

\noindent {\bf Theorem 2.} Let $0\le X\sim F.$  Then  the
following statements are true.\\
(a) If $s\in (2n-1, 2n+1)$ for some integer $n\ge 1,$ the moment
\begin{eqnarray}
m_{s}={\bf E}[X^s]=\frac{2}{\pi}\Gamma(s+1)\cos((s-2n)\pi/2)
\int_0^{\infty}\frac{{\bf E}[h_n(tX)]}{t^{s+1}}\,dt\ \ (\hbox{finite or infinite}).
\end{eqnarray}
(b) If, in addition, $m_s<\infty$ for some $s\in(2n-1,2n+1)$ with
integer  $n\ge 1,$
\begin{eqnarray}
m_{s}=\frac{2}{\pi}\Gamma(s+1)\cos((s-2n)\pi/2)\int_0^{\infty}\frac{H_n(t)}{t^{s+1}}\,dt.
\end{eqnarray}
(c) If   $m_s<\infty$ for some $s\in(0,1),$
\begin{eqnarray}
m_{s}=\frac{2}{\pi}\Gamma(s+1)\cos(s\pi/2)\int_0^{\infty}\frac{{\bf E}[\sin (tX)]}{t^{s+1}}\,dt
=\frac{2}{\pi}\Gamma(s+1)\cos(s\pi/2)\int_0^{\infty}\frac{H_0(t)}{t^{s+1}}\,dt.
\end{eqnarray}
(d) If $0<X\sim F,$ the improper integral
\begin{eqnarray}\int_0^{\infty}\frac{{\bf E}[\sin
(tX)]}{t}\,dt=\int_0^{\infty}\frac{\hbox{Im}(\phi(t))}{t}\,dt=\frac{\pi}{2},
\end{eqnarray}
regardless of the distribution $F.$\vspace{0.3cm}

 In view of (14), we can define
a density function
\begin{eqnarray}g(t)=\frac{2}{\pi}\,\frac{{\bf E}[\sin (tX)]}{t},\ \ \ t>0,
\end{eqnarray} provided
${\bf E}[\sin (tX)]\ge 0$ on $(0,\infty).$ This leads to the
following reciprocal moment relation. Two illustrative examples are
given in Remark 3 below. See also Remark 4 and Laue (1986) for other
moment relations by using the mixture of exponential distribution.

\noindent {\bf Theorem 3.} Let $0<X\sim F$ with ${\bf E}[\sin
(tX)]\ge 0$ for all $t\in{\mathbb R}_+,$ and let $0<Y\sim G$ with
density function $g$ defined in (15).\\ (a) If ${\bf E}[X^s]<\infty$
for some $s\in(0, 1),$ then
\begin{eqnarray}
{\bf E}[Y^{-s}]=\frac{1}{\Gamma(s+1)\cos(s\pi/2)}\,{\bf E}[X^s].
\end{eqnarray}
(b) If, instead, ${\bf E}[X^s]<\infty$ for all $s\in(-1, 1),$ then
(16) holds true for all $s\in(-1, 1).\vspace{0.3cm}$

It is surprising that the condition ${\bf E}[\sin (tX)]\ge 0$ in
Theorem 3 is fulfilled by many distributions as shown in Lemma 3 and
Corollary 1  below.

 \noindent {\bf Corollary 1.} Let $Z$ have a P\'olya-type ch.f.\,\,$\varphi$
 and let $X>0$ obey the distribution $F=1-\varphi.$ Assume further that
${\bf E}[X^s]<\infty$ for all $s\in(-1, 1).$ Then
\begin{eqnarray*}
{\bf E}[|Z|^{-s}]=\frac{1}{\Gamma(s+1)\cos(s\pi/2)}\,{\bf E}[X^s],\ \ |s|<1.
\end{eqnarray*}

Similarly, for length-biased distributions, we have the following
reciprocal moment relation by Theorem 1$^{\prime}.$ \\
 \noindent {\bf Theorem 4.} Let $0\le X\sim F$ with mean $m_1\in(0,\infty).$
Define the random variable $Z\sim H$ with the length-biased
distribution $H(z)=\frac{1}{m_1}\int_0^zt\,dF(t),\,z\ge 0,$ and the
random variable $Y\sim G$ with the density
$$g(y)=\frac{2}{\pi}\,\frac{1-{\bf E}[\cos(yX)]}{m_1y^2},\ \ y>0.$$
If ${\bf E}[X^s]<\infty$ for some $s\in(0,2),$ then
\begin{eqnarray*}
{\bf E}[Y^{1-s}]=\frac{1}{\Gamma(s+1)\sin(s\pi/2)}\,\frac{{\bf E}[X^{s}]}{m_1}
=\frac{1}{\Gamma(s+1)\sin(s\pi/2)}\,{\bf E}[Z^{s-1}].\vspace{0.5cm}
\end{eqnarray*}

 For nonnegative random variables $X\sim F$, we define
 the Laplace--Stieltjes transform $L(\lambda)={\bf
E}[\exp(-\lambda X)],\ \lambda\ge 0.$  It is known that the
distribution $F$ of $X$ is uniquely determined by $L;$ more
precisely, for each continuity point $x\ge 0$ of $F,$
$$F(x)=\lim_{n\to\infty}\sum_{k\le nx}(-1)^k\frac{n^k}{k!}L^{(k)}(n)$$
(see also Remark 2 below). We can apply $L$ (instead of
Fourier--Stieltjes transform) to get the formulas for fractional
moments. To do this, let us define the auxiliary functions
\begin{eqnarray}
q_n(\lambda)=(-1)^{n+1}\biggl(e^{-\lambda}-\sum_{k=0}^n\frac{(-1)^k}{k!}\lambda^k\biggr)\ge
0,\ \lambda\ge 0,\ \ n\ge 0,
\end{eqnarray}
and
\begin{eqnarray} Q_n(\lambda)=(-1)^{n+1}\biggl(L({\lambda})-\sum_{k=0}^n\frac{(-1)^k}{k!}m_k\lambda^k\biggr)\ge
0,\ \lambda\ge 0,\ \ n\ge 0,
\end{eqnarray}
provided the moments are finite (see, e.g., Hu and Lin 2008).  Note
that the derivative $q_n^{\prime}(\lambda)=q_{n-1}(\lambda)$ for
$n\ge 1.$ Then we have the following result. See also Urbanik (1993,
p.\,326) and Klar (2003) for the formula (20) {\it under the
assumption of finite fractional moments}.

\noindent {\bf Theorem 5.} Let $0\le X\sim F$ and $n<s<n+1,$ where
$n\ge 0$ is an integer. Then the fractional moment (finite or
infinite)
\begin{eqnarray}
m_{s}=\frac{1}{\pi}\Gamma(s+1)\sin((s-n)\pi)\int_0^{\infty}\frac{{\bf E}[q_n({\lambda}X)]}{{\lambda}^{s+1}}\,d{\lambda}
=\frac{(-1)^{n+1}}{\Gamma(-s)}\int_0^{\infty}\frac{{\bf E}[q_n({\lambda}X)]}{{\lambda}^{s+1}}\,d{\lambda}.
\end{eqnarray}
If, in addition, $m_s<\infty$ for some $s\in(n, n+1)$ with integer
$n\ge 0,$ then
\begin{eqnarray}
m_{s}=\frac{1}{\pi}\Gamma(s+1)\sin((s-n)\pi)\int_0^{\infty}\frac{Q_n({\lambda})}{{\lambda}^{s+1}}\,d{\lambda}
=\frac{(-1)^{n+1}}{\Gamma(-s)}
\int_0^{\infty}\frac{Q_n({\lambda})}{{\lambda}^{s+1}}\,d{\lambda}.\vspace{0.3cm}
\end{eqnarray}

As shown in Lemma 5, the function $Q_n$ in (18) is  closely related
to the high-order equilibrium distribution  $F_{(n+1)}$ with respect
to (w.r.t.)  $F$ defined below. For $0\le X\sim F$ with mean
$m_1\in(0,\infty),$ we define the first-order equilibrium
distribution  by $F_{(1)}(x)=m_1^{-1}\int_0^x\overline{F}(y)\,dy,\
x\ge 0,$   where $\overline{F}=1-F.$ The high-order equilibrium
distributions are defined iteratively. Namely, the $n$th-order
equilibrium distribution (w.r.t.\,\,$F$) is
$F_{(n)}(x)=m_{(n-1)}^{-1}\int_0^x\overline{F}_{(n-1)}(y)\,dy,\ x\ge
0,$ provided the mean $m_{(n-1)}$ of ${F}_{(n-1)}$ is finite. It is
known that $m_{(n-1)}=m_n/(n\,m_{n-1})$ (see, e.g., Lin 1998b,
p.\,265, or Harkness and Shantaram 1969).
  An
application of Theorem 5 leads to the next result similar to
Theorems 3 and 4.

\noindent {\bf Theorem 6.} Let $0\le X\sim F$ with $n$th moment
$m_n\in(0,\infty)$ for some integer $n\ge 1.$  Let $X_{(n)}$ have
the $n$th-order equilibrium distribution $F_{(n)}$ and
Laplace--Stieltjes transform $L_{(n)}(\lambda)\equiv{\bf
E}[\exp(-\lambda X_{(n)})],\ \lambda\ge 0.$ Assume further that  $Y>0$ has distribution $G=1-L_{(n)}.$ \\
(a) If ${\bf E}[X^s]<\infty$ for some $s\in(n-1,n),$ then
\begin{eqnarray}
{\bf E}[Y^{n-s}]=\frac{n!(n-s)\pi}{m_n\Gamma(s+1)\sin((n-s)\pi)}\,{\bf E}[X^s].
\end{eqnarray}
(b) If, instead, ${\bf E}[X^s]<\infty$ for all $s\in(n-1, n+1),$
then (21) holds true for all $s\in(n-1, n+1).$ The constant is
defined for $s=n$ by continuity to be equal to $1/m_n.$

Finally, we present the counterpart result of Hall (1983) for
distributions on the right-half line, which claims that
 a {\it proper} subset
 of the previous translated moments in (6) is enough to characterize such a distribution.
 Theorem 7 treats the case of finite $s$th fractional moment with $s\in(0,1),$
 while Theorem 8 deals with the remaining cases which are more involved.

\noindent {\bf Theorem 7.}  Let $0\le X\sim F$ and let the
fractional moment ${\bf E}[X^s]\in(0,\infty)$ for some $s\in(0,1).$
Then the distribution  $F$ is completely determined by the sequence
$$M_s^+=\{{\bf E}[X^s]\} \cup\{{\bf
E}[(X+a_k)^s]\}_{k=1}^{\infty}$$  of translated fractional moments,
where $\{a_k\}_{k=1}^{\infty}$ is a sequence of positive and
distinct real
numbers satisfying one of the following conditions:\\
(a) $\lim_{k\to\infty}a_k=\infty$ and
$\sum_{k=1}^{\infty}1/a_k=\infty;$\\
(b) $\lim_{k\to\infty}a_k=a_0\in(0,\infty);$\\
(c) $\lim_{k\to\infty}a_k=0$ and $\sum_{k=1}^{\infty}a_k=\infty.$
\medskip\\
\noindent {\bf Theorem 8.}  Let $0\le X\sim F$ and let the
fractional moment ${\bf E}[X^s]\in(0,\infty)$ for some
$s\in(n,n+1),$ where $n\ge 1$ is an integer. Further, suppose that
$\{a_k\}_{k=1}^{\infty}$ is a sequence of positive and distinct real
numbers satisfying: $\lim_{k\to\infty}a_k=\infty$ and
$\sum_{k=1}^{\infty}1/a_k=\infty.$ Then the distribution  $F$ is
completely determined by the sequence $$M_s^{++}=\{{\bf E}[X^s]\}
\cup\{{\bf E}[(X+a_k)^s]\}_{k=1}^{\infty} \cup\{{\bf
E}[(X+2a_k)^s]\}_{k=1}^{\infty}$$  of translated fractional moments.

An immediate consequence of the last two theorems is the following.

 \noindent {\bf Corollary 2.} Let $0\le X\sim F.$\\
(a) If the fractional moment ${\bf E}[X^s]\in (0,\infty)$ for some
 $s\in(0,1),$ then the distribution $F$ is characterized  by any one of the two sequences:
\[(\hbox{i})\ \{{\bf E}[X^s]\} \cup\{{\bf
E}[(X+p_j)^s]\}_{j=1}^{\infty}\ \ \ \hbox{and}\ \ \ \hbox{(ii)}\
\{{\bf E}[X^s]\} \cup\{{\bf E}[(X+1/p_j)^s]\}_{j=1}^{\infty},\] where
$p_j$ is the $j$th prime
number.\\
(b) If, instead, the fractional moment ${\bf E}[X^s]\in (0,\infty)$
for some non-integer
 $s>1,$ then the distribution $F$ is characterized by  the sequence
\[\{{\bf E}[X^s]\} \cup\{{\bf E}[(X+p_j)^s]\}_{j=1}^{\infty}
\cup\{{\bf E}[(X+2p_j)^s]\}_{j=1}^{\infty}.\]

\noindent{\bf 3. Lemmas}

To prove the main results,  we need some preliminary lemmas in the
sequel.

\noindent {\bf Lemma 1.} Let  $n\ge 0$ be an integer and
$2n<s<2n+2.$ Then the integral
$$\int_0^{\infty}\frac{g_n(t)}{t^{s+1}}\,dt=\frac{\pi}{2}[\Gamma(s+1)\sin((s-2n)\pi/2)]^{-1}.$$
\noindent{\bf Proof.} Note that $g_0(t)=1-\cos t,\ \ t\in{\mathbb
R},$ and the derivatives $g_k^{\prime}(t)=h_k(t),\
h_k^{\prime}(t)=g_{k-1}(t)$ for $k\ge 1.$ Moreover, by L'Hospital
rule, we have, for $2\le 2k<r<2(k+1),$ the limits
$$ \lim_{t\to 0^+}\frac{g_k(t)}{t^r}=0=\lim_{t\to \infty}\frac{g_k(t)}{t^r},\ \  \ \
\lim_{t\to 0^+}\frac{h_k(t)}{t^{r-1}}=0=\lim_{t\to \infty}\frac{h_k(t)}{t^{r-1}}.$$
First, we recall the identity:
$$\int_{-\infty}^{\infty}\frac{1-\cos t}{|t|^{r+1}}\,dt=\frac{\pi}{\Gamma(r+1)\sin(r\pi/2)},\ \ r\in(0,2)$$
(see, e.g., Chung 2001, p.\,159), from which the case $n=0$ in the
lemma follows. Then, by integration by parts and induction, we have,
for $n\ge 1,$ the integral
\begin{eqnarray*}
\int_0^{\infty}\frac{g_n(t)}{t^{s+1}}\,dt&=&\frac{1}{s}\int_0^{\infty}
\frac{h_n(t)}{t^{s}}\,dt=\frac{1}{s(s-1)}\int_0^{\infty}\frac{g_{n-1}(t)}{t^{s-1}}\,dt\\
&=&\cdots\cdots\cdots\\
&=&\frac{1}{s(s-1)\cdots(s-2n+1)}\int_0^{\infty}\frac{g_{0}(t)}{t^{s-2n+1}}\,dt\\
&=&\frac{1}{s(s-1)\cdots(s-2n+1)}\int_0^{\infty}\frac{1-\cos t}{t^{s-2n+1}}\,dt\\
&=&\frac{\pi}{2\Gamma(s+1)\sin((s-2n)\pi/2)}.
\end{eqnarray*}
The proof is complete.

 \noindent {\bf Lemma 2.}
Let $n\ge 0$ be an integer and $2n-1<s<2n+1.$ Then the integral
$$\int_0^{\infty}\frac{h_n(t)}{t^{s+1}}\,dt=\frac{\pi}{2}[\Gamma(s+1)\cos((s-2n)\pi/2)]^{-1}.$$

Note that the above identity holds true in the sense of improper
integral (not the Lebesgue integral on $(0,\infty)$) when $n=0$ and $s\in(-1,0].$\\
 \noindent{\bf Proof of Lemma 2.} First,
we consider the case $n\ge 1$ and $2n-1<s<2n+1$, or, equivalently,
$2\le 2n<s+1<2n+2.$ Then from the proof of Lemma 1 (replacing $s$ by
$s+1$), we have the integral
$$\int_0^{\infty}\frac{h_n(t)}{t^{s+1}}\,dt=\frac{(s+1)\pi}{2\Gamma(s+2)\sin((s+1-2n)\pi/2)}
=\frac{\pi}{2\Gamma(s+1)\cos((s-2n)\pi/2)}.$$
Therefore, it remains to prove the case $n=0.$ Namely, for $-1<s<1,$
we want to prove
$$\int_0^{\infty}\frac{h_0(t)}{t^{s+1}}\,dt=\int_0^{\infty}\frac{\sin t}{t^{s+1}}\,dt
=\frac{\pi}{2\Gamma(s+1)\cos(s\pi/2)}.$$ The case $s=0$ is trivial.
For $s\ne 0,$ we recall the identity:
$$\int_0^{\infty}\frac{\sin t}{t^{s+1}}\,dt=\Gamma(-s)\sin(-s\pi/2),\ \ 0<|s|<1$$
(see, e.g., Gradshteyn and Ryzhik 2014, Formulas 3.721(1) and
3.761(4)). Finally, Euler's reflection formula,
$\Gamma(1-z)\Gamma(z)={\pi}/{\sin(z\pi)}$ for non-integer $z,$
completes the proof.

 \noindent {\bf Lemma 3.} Let $0<X\sim F$ with density function $f(x)$ decreasing to zero as $x\to\infty.$
 Then ${\bf E}[\sin (tX)]\ge 0$ for all $t\in{\mathbb R}_+.$\\
{\bf Proof.} See, e.g., Lukacs (1970), p.\,84.\\
 \noindent {\bf Lemma 4.} Let  $n\ge 0$ be an integer and $n<s<n+1.$
Then the integral
$$\int_0^{\infty}\frac{q_n(\lambda)}{\lambda^{s+1}}\,d\lambda=
\frac{\Gamma(n+1-s)}{s(s-1)(s-2)\cdots(s-n)}=(-1)^{n+1}\Gamma(-s).$$
\noindent{\bf Proof.} Note that $q_0(\lambda)=1-\exp(-\lambda), \
\lambda\ge 0,$ and that the derivative
$q_k^{\prime}(\lambda)=q_{k-1}(\lambda)$ for $k\ge 1.$ Moreover, by
L'Hospital rule, we have, for $0\le k<r<k+1,$ the limits
$$ \lim_{\lambda\to 0^+}\frac{q_k(\lambda)}{\lambda^r}=0=\lim_{\lambda\to \infty}\frac{q_k(\lambda)}{\lambda^r}.$$
First, it can be shown that for $0<r<1,$ we have the integral
$$\int_0^{\infty}\frac{q_0(\lambda)}{\lambda^{r+1}}\,d\lambda=
\int_0^{\infty}\frac{1-e^{-\lambda}}{\lambda^{r+1}}\,d\lambda=\frac{\Gamma(1-r)}{r}.$$
This is exactly the case $n=0.$ Then, by integration by parts and
induction, we have, for $n\ge 1,$ the integral
\begin{eqnarray*}
\int_0^{\infty}\frac{q_n(\lambda)}{\lambda^{s+1}}\,d\lambda
&=&\frac{1}{s}\int_0^{\infty}\frac{q_{n-1}(\lambda)}{\lambda^{s}}\,d\lambda
=\ \cdots\cdots\cdots\\
&=&\frac{1}{s(s-1)\cdots(s-n+1)}\int_0^{\infty}\frac{q_{0}(\lambda)}{\lambda^{s-n+1}}\,d\lambda\\
&=&\frac{\Gamma(n+1-s)}{s(s-1)\cdots(s-n)}=(-1)^{n+1}\Gamma(-s).
\end{eqnarray*}
The proof is complete.

The next lemma is an application of Lin's (1994) theorem to the
specific function $g(x)=\exp(-\lambda x),\ x\ge 0,$ where $\lambda >
0$ is a constant.

\noindent {\bf Lemma 5.} Let $0\le X\sim F$ with $n$th moment
$m_n\in(0,\infty)$  for some integer $n\ge 1,$ and let $X_{(n)}$
have the $n$th-order equilibrium distribution $F_{(n)}.$
 Then the Laplace--Stieltjes transform of $X_{(n)}$
is equal to
$$L_{(n)}(\lambda)={\bf E}[\exp(-\lambda X_{(n)})]=n!Q_{n-1}(\lambda)/(m_n\lambda^n),\ \ \lambda\ge 0,$$
where the function $Q_{n-1}$ is defined in (18), and the RHS of the
above equality is defined for $\lambda=0$ by continuity to be equal
to $1.$

For a proof of the following variant of the M\"untz--Sz\'asz
Theorem,
 see Lin (1993), where  $`\stackrel{d}{=}$ '
means `equal in distribution'.

 \noindent {\bf Lemma 6.} Let $X_1$
and $X_2$ be two nonnegative random variables with
Laplace--Stieltjes transforms $L_1$ and $L_2$, respectively. Let
$\{\lambda_k\}_{k=1}^{\infty}$ be a sequence of positive and
distinct real
numbers satisfying one of the following conditions:\\
(a) $\lim_{k\to\infty}\lambda_k=\infty$ and
$\sum_{k=1}^{\infty}1/\lambda_k=\infty;$\\
(b) $\lim_{k\to\infty}\lambda_k=\lambda_0\in(0,\infty);$\\
(c) $\lim_{k\to\infty}\lambda_k=0$ and
$\sum_{k=1}^{\infty}\lambda_k=\infty.$\\
Further, assume that $L_1(\lambda_k)=L_2(\lambda_k)$ for all $k\ge
1.$ Then $X_1\stackrel{d}{=}X_2.$

\newpage
\noindent{\bf 4.  Proofs of main results}

\noindent{\bf Proof of Theorem 1$^{\prime}$.} It suffices to prove
formula (9). Since the even function $g_n\ge 0$ for all $n\ge 0,$
Tonelli's Theorem applies (see, e.g., Royden 1988, p.\,309) and we
have, by $g_n(0)=0$ and Lemma 1,
\begin{eqnarray*}
\int_0^{\infty}\frac{{\bf E}[g_n(tX)]}{t^{s+1}}\,dt
&=&\frac{1}{2}\int_{-\infty}^{\infty}\frac{{\bf E}[g_n(tX)]}{|t|^{s+1}}\,dt
=\frac{1}{2}\int_{-\infty}^{\infty}\int_{-\infty}^{\infty}\frac{g_n(tx)}{|t|^{s+1}}dF(x)\,dt\\
&=&\frac{1}{2}\int_{-\infty}^{\infty}\int_{-\infty}^{\infty}\frac{g_n(tx)}{|t|^{s+1}}dt\,dF(x)
=\frac{1}{2}\int_{-\infty}^{\infty}|x|^s\int_{-\infty}^{\infty}\frac{g_n(t)}{|t|^{s+1}}dt\,dF(x)\\
&=&\int_{0}^{\infty}\frac{g_n(t)}{t^{s+1}}dt\int_{-\infty}^{\infty}|x|^sdF(x)
=\frac{\pi}{2}[\Gamma(s+1)\sin((s-2n)\pi/2)]^{-1}\mu_s.
\end{eqnarray*}
This proves (9), and hence (10) follows. The
proof is complete.\medskip\\
\noindent{\bf Proof of Theorem 2.} We first prove formula (11).
Since the odd function $h_n(t)\ge 0$ on ${\mathbb R}_+$ for $n\ge
1,$ Tonelli's Theorem applies and we have, by $h_n(0)=0$ and Lemma
2,
\begin{eqnarray*}
\int_0^{\infty}\frac{{\bf E}[h_n(tX)]}{t^{s+1}}\,dt
&=&\int_{0}^{\infty}\int_{0}^{\infty}\frac{h_n(tx)}{t^{s+1}}dF(x)\,dt
=\int_{0}^{\infty}\int_{0}^{\infty}\frac{h_n(tx)}{t^{s+1}}dt\,dF(x)\\
&=&\int_{0}^{\infty}x^s\int_{0}^{\infty}\frac{h_n(t)}{t^{s+1}}dt\,dF(x)
=\int_{0}^{\infty}\frac{h_n(t)}{t^{s+1}}dt\int_{0}^{\infty}x^sdF(x)\\
&=&\frac{\pi}{2}[\Gamma(s+1)\cos((s-2n)\pi/2)]^{-1}m_s.
\end{eqnarray*}
This proves part (a), and hence part (b) with $n\ge 1$ follows.

As for part (c), namely, the case $n=0$ and $s\in(0, 1),$ we first
prove that
\begin{eqnarray}\int_0^{\infty}\frac{|h_0(t)|}{t^{s+1}}\,dt=\int_0^{\infty}\frac{|\sin
t|}{t^{s+1}}\,dt=\int_0^{1}\frac{|\sin
t|}{t^{s+1}}\,dt+\int_1^{\infty}\frac{|\sin
t|}{t^{s+1}}\,dt\equiv I_1+I_2<\infty.
\end{eqnarray}
Since $\lim_{t\to 0+}(\sin t)/t=1,\ \sup_{t\in(0,1]}|\sin t|/t\equiv
M<\infty,$ and hence the first integral $I_1\le
M\int_0^{1}t^{-s}dt={M}/{(1-s)}<\infty.$ On the other hand, the
second  integral
$$I_2\le \int_1^{\infty}\frac{1}{t^{s+1}}\,dt=1/s<\infty,$$ which completes the proof of (22).
This in turn implies that
$$\int_0^{\infty}\int_0^{\infty}\frac{|\sin(tx)|}{t^{s+1}}dt\,dF(x)=
\int_0^{\infty}x^s\int_0^{\infty}\frac{|\sin t|}{t^{s+1}}dt\,dF(x)={\bf E}[X^s]\int_0^{\infty}\frac{|\sin
t|}{t^{s+1}}\,dt<\infty,$$ provided $m_s<\infty.$
 Hence, the Fubini--Tonelli Theorem
applies (see, e.g., Rudin 1987, p.\,164) and (13) holds true by
Lemma 2.  Part (d) is an immediate consequence of (2) or (4) by
letting $x=0.$ The proof is complete.
\medskip\\
\noindent{\bf Proof of Theorem 3.}  To prove part (a), note that, by
Theorem 2(c), the moment
$${\bf E}[Y^{-s}]=\int_0^{\infty}t^{-s}g(t)\,dt=\int_0^{\infty}
\frac{2}{\pi}\,\frac{{\bf E}[\sin (tX)]}{t^{s+1}}\,dt
=\frac{1}{\Gamma(s+1)\cos(s\pi/2)}\,{\bf E}[X^s].
$$
Part (b) follows from part (a) with the help of analytic
continuation of Mellin transform.
\medskip\\
\noindent{\bf Proof of Corollary 1.} Since $\varphi$ is a
P\'olya-type ch.f., $Z$ has a P\'olya-type density function $p$,
which is symmetric  on ${\mathbb R}$ and of the form
$$p(z)=\frac{1}{\pi}\int_0^{\infty}\cos(zt)\varphi(t)\,dt=\frac{1}{\pi}\int_0^{\infty}\cos(zt)[1-F(t)]\,dt,\ \
z>0$$ (Lukacs 1970, pp.\,\,83--84). By integration by parts, we rewrite the density
$$p(z)=\frac{1}{\pi z}\int_0^{\infty}\sin(zt)\,dF(t)=\frac{1}{\pi z}{\bf E}(\sin(zX)\ge 0,\ \ z>0.$$
Note that  $Y=|Z|$ has the density function $g(y)=2p(y),\ y>0.$
Therefore, the required result follows from Theorem 3. The proof is
complete.

\noindent{\bf Proof of Theorem 4.} Note that, by Theorem
1$^{\prime}$ with $n=0,$ the moment
$${\bf E}[Y^{1-s}]=\int_0^{\infty}y^{1-s}g(y)\,dy=\int_0^{\infty}
\frac{2}{\pi}\,\frac{1-{\bf E}[\cos(yX)]}{m_1y^{s+1}}\,dy
=\frac{1}{\Gamma(s+1)\sin(s\pi/2)}\,\frac{{\bf E}[X^{s}]}{m_1}.$$

\noindent{\bf Proof of Theorem 5.} It suffices to prove formula
(19). Since the function $q_n\ge 0,$ Tonelli's Theorem applies and
we have, by $q_n(0)=0$ and Lemma 4,
\begin{eqnarray*}
\int_0^{\infty}\frac{{\bf E}[q_n(\lambda X)]}{\lambda^{s+1}}\,d\lambda
&=&\int_{0}^{\infty}\int_{0}^{\infty}\frac{q_n(\lambda x)}{\lambda^{s+1}}dF(x)\,d\lambda
=\int_{0}^{\infty}\int_{0}^{\infty}\frac{q_n(\lambda x)}{\lambda^{s+1}}d\lambda\, dF(x)\\
&=&\int_{0}^{\infty}x^s\int_{0}^{\infty}\frac{q_n(\lambda )}{\lambda^{s+1}}d\lambda\, dF(x)
=\int_{0}^{\infty}\frac{q_n(\lambda )}{\lambda^{s+1}}d\lambda\int_{0}^{\infty}x^sdF(x)\\
&=&(-1)^{n+1}\Gamma(-s)\,m_s.
\end{eqnarray*}
Therefore, (19) holds true, and (20) follows. The proof is complete.

 \noindent{\bf Proof of Theorem
6.} For $s\in(n-1,n),$ let us consider the moment
\begin{eqnarray*}{\bf
E}[Y^{n-s}]&=&\int_0^{\infty}\Pr(Y^{n-s}>y)\,dy=\int_0^{\infty}\Pr(Y>y^{1/(n-s)})\,dy \\
&=&\int_0^{\infty}\overline{G}(y^{1/(n-s)})\,dy=
{(n-s)}\int_0^{\infty}\frac{\overline{G}(\lambda)}{\lambda^{s-n+1}}\,d\lambda\\
&=&{(n-s)}\int_0^{\infty}\frac{L_{(n)}(\lambda)}{\lambda^{s-n+1}}\,d\lambda=
\frac{n!(n-s)}{m_n}\int_0^{\infty}\frac{Q_{n-1}(\lambda)}{\lambda^{s+1}}\,d\lambda.
\end{eqnarray*} The last equality is due to Lemma 5.
Therefore, part (a) follows from Theorem 5, while part (b) follows
from part (a) with the help of analytic continuation of Mellin
transform.

 \noindent{\bf Proof of Theorem
7.} Suppose that the random variables $0\le X\sim F$ and $0\le Y\sim
G$ satisfy (i) ${\bf E}[X^s] ={\bf E}[Y^s]=m_s\in(0,\infty)$ for
some $s\in(0,1)$ and
\begin{eqnarray}(\hbox{ii})~~{\bf
E}[(X+a_k)^s]={\bf E}[(Y+a_k)^s],\ \ k=1,2,\ldots,
\end{eqnarray} where the
sequence $\{a_k\}_{k=1}^{\infty}$ is defined in the theorem. Then we
want to prove that $F=G.$

Denote the constant $C_s=\frac{1}{\pi}\Gamma(s+1)\sin(s\pi).$ Then,
by Theorem 5, we have  the moment
\begin{eqnarray}{\bf
E}[X^s]=C_s\int_0^{\infty}\frac{1-L_X(\lambda)}{\lambda^{s+1}}\,d\lambda,\end{eqnarray}
where $L_X$ is the Laplace--Stieltjes transform of $X\sim F.$
 In other words, the function
$$h_X(\lambda)=\frac{C_s}{m_s}\,\frac{1-L_X(\lambda)}{\lambda^{s+1}}\ge 0,\ \lambda> 0,$$
is a bona fide density. Let $Z_X\sim H_X$ have the density $h_X.$
Similarly, let $Z_Y\sim H_Y$ have the corresponding density (derived
from $Y\sim G$)
$$h_Y(\lambda)=\frac{C_s}{m_s}\,\frac{1-L_Y(\lambda)}{\lambda^{s+1}}\ge 0,\ \lambda> 0.$$
\indent Since $L_{X+a_k}(\lambda)={\bf
E}[\exp(-\lambda(X+a_k))]=\exp(-a_k\lambda)L_{X}(\lambda),\,\lambda\ge
0,$  it follows from (23) and (24)  that
$$\int_0^{\infty}e^{-a_k\lambda}\,\frac{L_Y(\lambda)-L_X(\lambda)}{\lambda^{s+1}}\,d\lambda=0,\ \ k=1,2,\ldots,$$
or, equivalently,
\begin{eqnarray*}\int_0^{\infty}e^{-a_k\lambda}\,dH_X(\lambda)
&=&\int_0^{\infty}e^{-a_k\lambda}h_X(\lambda)\,d\lambda=
\int_0^{\infty}e^{-a_k\lambda}h_Y(\lambda)\,d\lambda\\
&=&\int_0^{\infty}e^{-a_k\lambda}\,dH_Y(\lambda),\
\ k=1,2,\ldots.
\end{eqnarray*}
This implies that $H_X=H_Y$ by Lemma 6, and hence $h_X=h_Y,$ which
in turn implies that  $L_X=L_Y.$ Therefore, $X\stackrel{d}{=}Y$ by
the uniqueness theorem for Laplace--Stieltjes
transforms. The proof is complete.\medskip\\
\noindent{\bf Proof of Theorem 8.} Suppose that the random variables
$0\le X\sim F$ and $0\le Y\sim G$ satisfy (i) ${\bf E}[X^s] ={\bf
E}[Y^s]=m_s\in(0,\infty)$ for some $s\in(n,n+1),$ where $n\ge 1$ is
a fixed integer, and
\begin{eqnarray}(\hbox{ii})~~{\bf
E}[(X+a)^s]={\bf E}[(Y+a)^s],\ \ a=a_k,\ 2a_{k},\ \ k=1,2,\ldots.
\end{eqnarray}  Then we
want to prove that $F=G.$

Firstly, we show that
\begin{eqnarray}{\bf E}[X^{k}]={\bf E}[Y^{k}]=m_{k},\ \
1\le k\le n.\end{eqnarray}
 Mimicking the proof of Hall (1983), we have the relation: for
all $0\le\ell<n,$
\begin{eqnarray*}
& &a^{\ell+1}\{2a^{-s}{\bf
E}[(X+a)^s]-2^{\ell+1}(2a)^{-s}{\bf
E}[(X+2a)^s]\}\\
&=&\sum_{j=0}^{\ell}{s\choose j}a^{\ell+1-j}(2-2^{\ell+1-j}){\bf
E}[X^j]+{s\choose \ell+1}{\bf
E}[X^{\ell+1}] + {o}(1)\ \ \ {\hbox{as}}\ \ 0<a\to\infty.
\end{eqnarray*}
Therefore, the sequence $M_s^{++}$ determines the moments ${\bf
E}[X^{k}],\  1\le k\le n,$ and (26) holds true due to the assumption
(25).

Secondly, let $L_X$ and $L_Y$ be the Laplace--Stieltjes transforms
of $X$ and $Y,$ respectively. Denote the constant
$C_s=\frac{1}{\pi}\Gamma(s+1)\sin((s-n)\pi),$ and define the
function
$$h_X(\lambda)=\frac{C_s}{m_s}\,\frac{Q_n({\lambda})}{{\lambda}^{s+1}}\,\ge 0,\ \ \lambda>0,$$
in which $L=L_X$ in Theorem 5. Then it follows from (20) that
$\int_0^{\infty}h_X(\lambda)d{\lambda}=1,$ namely, $h_X$ is a bona
fide density function. Similarly,   define the density function
$$h_Y(\lambda)=\frac{C_s}{m_s}\,\frac{Q_n^*({\lambda})}{{\lambda}^{s+1}}\,\ge 0,\ \ \lambda>0,$$
where
$$Q_n^*(\lambda)=(-1)^{n+1}\biggl(L_Y({\lambda})-\sum_{k=0}^n\frac{(-1)^k}{k!}m_k\lambda^k\biggr)\ge
0,\ \lambda\ge 0. $$ \indent Thirdly, recall from (19) that for
$a>0,$ the translated moment \begin{eqnarray}{\bf
E}[(X+a)^s]&=&C_s\int_0^{\infty}\frac{{\bf
E}[q_n({\lambda}(X+a))]}{{\lambda}^{s+1}}\,d{\lambda}\nonumber\\
&=&
C_s(-1)^{n+1}\int_0^{\infty}\frac{1}{\lambda^{s+1}}\biggl(e^{-a\lambda}L_X(\lambda)-\sum_{k=0}^n\frac{(-1)^k}{k!}{\bf
E}[(X+a)^k]\biggr)\,d{\lambda}.
\end{eqnarray}
Combining (25) through (27) yields
$$\int_0^{\infty}e^{-a_k\lambda}\,\frac{L_X(\lambda)-L_Y(\lambda)}{\lambda^{s+1}}\,d\lambda =0,\ \ k=1,2,\ldots,$$
or, equivalently,
$$\int_0^{\infty}e^{-a_k\lambda}h_X(\lambda)\,d\lambda=
\int_0^{\infty}e^{-a_k\lambda}h_Y(\lambda)\,d\lambda,\ \ k=1,2,\ldots.$$
By Lemma 6, we finally conclude  that $h_X=h_Y,$ and hence
$L_X=L_Y,$
 $X\stackrel{d}{=}Y.$
The proof is complete.
\medskip\\
\noindent{\bf 5. Remarks}\\
\noindent{\bf Remark 1.} Note that the expression (5) improves
R\'enyi's (1970, p.\,369) result by deleting the redundant condition
$\mu_2<\infty.$ Hsu (1951) also proved the curious result that the
absolute moment (of negative order) $\mu_{-1}=\infty,$ if the
ch.f.\,\,$\phi$ has a nonnegative real part. The latter condition is
clearly satisfied by all the P\'olya-type ch.f.s (Lukacs 1970,
p.\,83).

\noindent{\bf Remark 2.} Suppose that  $0<X\sim F$ has a continuous
density $f$ and that the power function $f^r$ is integrable for some
$r>1.$ Then $f$ has the following representation in terms of the
Laplace--Stieltjes transform $L:$
$$f(x)=\lim_{n\to\infty}\frac{(-1)^{n-1}}{(n-1)!}(n/x)^nL^{(n-1)}(n/x),\ \ x>0$$
(see, e.g., Chung 2001, p.\,140).

\noindent{\bf Remark 3.} To illustrate the use of Theorem 3, let
$X\sim F$ have the standard exponential distribution
$F(x)=1-\exp(-x),\ x>0.$ Then $Y\sim G$ has the half-Cauchy density
$g(t)=(2/\pi)(1+t^2)^{-1},\,t>0,$ and hence, by (16), the moment
relation
$${\bf
E}[Y^{-s}]=\frac{1}{\Gamma(s+1)\cos(s\pi/2)}\,{\bf
E}[X^s]=\frac{1}{\cos(s\pi/2)},\,\ |s|<1.$$ More generally, let
$X\sim F$ have the Gamma density
$$f(x)=\frac{1}{\Gamma(\alpha)\beta^{\alpha}}\,x^{\alpha-1}\exp(-x/\beta),\ \ x>0,$$
where $\alpha\in(0,2]$ and $\beta>0.$ Then ${\bf E}[\sin (tX)]\ge
0,\ t\ge 0,$ and $Y\sim G$ has the density
$$g(t)=\frac{2}{\pi}\,\frac{{\bf E}[\sin (tX)]}{t}=\frac{2}{\pi}\,\frac{\sin(\alpha\theta)}{t(1+\beta^2t^2)^{\alpha/2}},\,\ t>0,$$
where $\theta=\arctan(\beta t)\in(0,\pi/2).$ Thus, we have the
moment relation
$${\bf
E}[Y^{-s}]=\frac{1}{\Gamma(s+1)\cos(s\pi/2)}\,{\bf
E}[X^s]=\frac{1}{\Gamma(s+1)\cos(s\pi/2)}\frac{\Gamma(s+\alpha)}{\Gamma(\alpha)}\beta^s,\,\
\ |s|<1.$$ \noindent{\bf Remark 4.} Let $0< X\sim F$ with Laplace--Stieltjes transform $L,$ and define the
random variable $Z\sim H$ with the mixture of exponential
distribution $H(\lambda)=\int_0^{\infty}(1-\exp(-\lambda
x))dF(x)=1-L(\lambda),\,\ \lambda\ge 0.$ Then we can improve Laue's
(1986, Lemma 3.2) moment relation as follows:
$${\bf E}[Z^{-s}]={\bf E}[{\cal E}^{-s}]{\bf E}[X^{s}]
\ \ \ (\hbox{finite or infinite}),\ \ \ s\in{\mathbb R},$$ where the
random variable ${\cal E},$ independent of $X,$ obeys the standard
exponential distribution and the moment ${\bf E}[{\cal E}^{-s}]$
equals $\Gamma(1-s)$ or $\infty$ according as $s<1$ or $s\ge 1.$ To
see this, recall that $Z\stackrel{d}{=}{\cal E}/X$ by the definition
of $H$ and hence
$${\bf E}[Z^{-s}]={\bf E}[({\cal E}/X)^{-s}]={\bf
E}[{\cal E}^{-s}X^s]={\bf E}[{\cal E}^{-s}]{\bf E}[X^s],\ \
\ s\in{\mathbb R},$$ due to the Monotone Convergence Theorem (see
Lin 1998a, Lemma 4).

\noindent{\bf Remark 5.} It seems hard to present the absolute
moments of {\it negative} order in terms of ch.f.\,(see, e.g., Hsu
1951, (15)). However, for {\it positive} random variable $X\sim F$
with Laplace--Stieltjes transform $L,$ we have the elegant formula\
\,
\begin{eqnarray}{\bf
E}[X^{-s}]=({\Gamma(s)})^{-1}\int_0^{\infty}\lambda^{s-1}L(\lambda)\,d\lambda,\
 s>0,
 \end{eqnarray} by the Tonelli Theorem (see, e.g.,  Pitman and Yor 2003,
 p.\,298, or Lin 1998b, p.\,264). Moreover, suppose that
$L(\lambda)$ is of the form
$L(\lambda)=\sum_{k=0}^{\infty}\frac{\varphi(k)}{k!}(-\lambda)^k$
for some suitable analytic function $\varphi(z),$ defined on the
half-plane ${\mathbb H}(\delta)=\{z\in{\mathbb C}:\hbox{Re}\ z\ge
-\delta\},$ where $\delta\in(0,1).$ Then when $s\in(0,\delta),$ the
formula (28) reduces to the following:
\[
{\bf
E}[X^{-s}]=({\Gamma(s)})^{-1}\int_0^{\infty}\lambda^{s-1}L(\lambda)\,d\lambda
=({\Gamma(s)})^{-1}[{\Gamma(s)}\varphi(-s)]=\varphi(-s),
\]
in which the penultimate equality follows from Ramanujan's Master
Theorem (see, e.g., Amdeberhan et al.\,2012 and Hardy 1999).

 \noindent{\bf
Acknowledgments.} The authors would like to thank Professor Lutz
Mattner (Universit\"at Trier, Germany) for mentioning of the related
references Riesz (1938), Rudin (1976) and Mattner (1992). The
comments and suggestions of Professor Jordan Stoyanov (Newcastle
University, UK) are also appreciated.

\vspace{0.5cm}\centerline{\bf References}

\begin{description}

\item Amdeberhan, T., Espinosa, O., Gonzalez, I., Harrison, M., Moll,
V.\,H. and Straub, A. (2012). Ramanujan's master theorem. {\it
Ramanujan J.}, {\bf 29}, 103--120.

\item Barndorff-Nielsen, O.\,E. and Stelzer, R. (2005). Absolute
moments of generalized hyperbolic distributions and approximate
inverse Gaussian L\'evy processes. {\it Scand. J. Statist.},
{\bf 32}, 617--637.

\item Brown, B.\,M. (1970). Characteristic functions, moments, and the central limit theorem.
{\it Ann. Math. Statist.}, {\bf 41}, 658--664.

\item  Brown, B.\,M. (1972). Formulae for absolute moments. {\it J. Aust. Math. Soc.}, {\bf 13}, 104--106.

\item {Chow, Y.\,S. and  Teicher, H}. (1997). {\it Probability Theory: Independence,
Interchangeability, Martingales}, 3rd edn. Springer, New York.

\item Chung, K.\,L. (2001). {\it A Course  in Probability Theory}, 3rd edn. Academic Press,
San Diego.

\item Gradshteyn, I.\,S. and Ryzhik, I.\,M. (2014).
{\it Table of Integrals, Series, and Products}, 8th edn.  Academic Press, New York.

\item Hall, P. (1983). A distribution is completely determined by its translated moments.
{\it Z. Wahrsch. Verw. Gebiete}, {\bf 62}, 355--359.

\item Hardy, G.\,H. (1999). {\it Ramanujan. Twelve Lectures on Subjects Suggested by
His Life and Work}.   AMS Chelsea Pub., Providence, R.I.

\item Harkness, W.\,L. and Shantaram, R. (1969). Convergence of a
sequence of transformations of distribution functions. {\it
Pacific J. Math.}, {\bf 31}, 403--415.

\item Hsu, P.\,L. (1951). Absolute moments and characteristic functions.
{\it J. Chinese Math. Soc. (New Series)}, {\bf 1}, 257--280.
Also in: {\it Pao-Lu Hsu Collected Papers}, edited by Kai Lai
Chung, pp.\,315--329. Springer,  New York, 1983.

\item Hu, C.-Y. and Lin, G.\,D. (2005). Some inequalities for
characteristic functions. {\it J. Math. Anal. Appl.}, {\bf 309},
336--352.

\item Hu, C.-Y. and Lin, G.\,D. (2008). Some inequalities for Laplace
transforms. {\it J. Math. Anal. Appl.}, {\bf 340}, 675--686.

\item  Kawata, T. (1972). {\it Fourier Analysis in Probability
Theory.} Academic Press,  New York.

\item Klar, B. (2003). On a test for exponentiality against Laplace order dominance.
{\it Statistics}, {\bf 37}, 505--515.

\item  Laue, G. (1983). Existence and representation of density
functions. {\it Math. Nachr.}, {\bf 114}, 7--21.

\item  Laue, G. (1986). Results on moments of non-negative random
variables. {\it Sankhy$\overline{a}$ Ser.\,A}, {\bf 48},
299--314.

\item Lin, G.\,D. (1993). Characterizations of the exponential
distribution via the blocking time in a queueing system. {\it
Statist. Sinica}, {\bf 3}, 577--581.

\item Lin, G.\,D. (1994). On a probabilistic generalization of Taylor's
theorem. {\it Statist. Probab. Lett.}, {\bf 19}, 239--243.

\item Lin, G.\,D. (1998a). On the Mittag-Leffler distributions. {\it J.
Statist. Plann. Inference}, {\bf 74}, 1--9.

\item Lin, G.\,D. (1998b). Characterizations of the ${\cal L}$-class of life
distributions. {\it Statist. Probab. Lett.}, {\bf 40}, 259--266.

\item Lukacs, E. (1970). {\it Characteristic Functions,} 2nd edn. Griffin, London.

\item Matsui, M. and Pawlas, Z. (2016). Fractional absolute moments of heavy tailed distributions.
{\it Braz. J. Probab. Statist.}, {\bf 30}, 272--298.

\item Mattner, L. (1992). Completeness of location families, translated moments, and
uniqueness of charges. {\it Probab. Theory Related Fields}, {\bf
92}, 137--149.

\item Nabeya, S. (1951). Absolute moments in 2-dimensional normal distribution.
{\it Ann. Inst. Statist. Math.}, {\bf 3}, 2--6.

\item Nabeya, S. (1952). Absolute moments in 3-dimensional normal
distribution. {\it Ann. Inst. Statist. Math.}, {\bf 4}, 15--30.

\item Pitman, J. and Yor, M. (2003). Infinitely divisible laws
associated with hyperbolic functions. {\it Canad. J. Math.},
{\bf 55}, 292--330.

\item R\'enyi, A. (1970). {\it Probability Theory}. North-Holland, Amsterdam.

\item Riesz, M. (1938). Int\'egrales de Riemann--Liouville et potentiels.
{\it Acta Szeged Sect. Math.}, {\bf 9}, 1--42.

\item  Rossberg, H.-J., Jesiak, B. and Siegel, G. (1985). {\it
Analytic Methods of Probability Theory.} Akademie-Verlag,
Berlin.
\item Royden, H.\,L. (1988). {\it Real Analysis}, 3rd edn. Prentice Hall, New Jersey.

\item  Rudin, W. (1976). $L^p$-isometries and equimeasurability.
{\it Indiana Univ. Math. J.}, {\bf 25}, 215--228.

\item  Rudin, W. (1987). {\it Real and Complex Analysis},  3rd edn. McGraw-Hill, New
York.

\item Urbanik, K. (1993). Moments of sums of independent random variables. In:
{\it Stochastic Processes}. A festschrift in honour of Gopinath
Kallianpur (eds. S. Cambanis, J.\,K. Ghosh, R.\,L. Karandikar,
and P.\,K. Sen),  pp.\,321--328, Springer, New York, 1993.

\item von Bahr, B. (1965). On the convergence of moments in the central limit theorem.
{\it Ann. Math. Statist.}, {\bf 36}, 808--818.

\item von Bahr, B. and Esseen, C.-G. (1965). Inequalities for the $r$th absolute moment
of a sum of random variables, $1\le r \le 2.$ {\it Ann. Math.
Statist.}, {\bf 36},  299--303.

\item  Zolotarev, V.\,M. (1957). Mellin--Stieltjes transforms in
probability theory. {\it Theory  Probab. Appl.}, {\bf 2},
433--460.

\end{description}

\end{document}